\documentclass[12pt,leqno,a4paper]{amsart}
\usepackage{amssymb,enumerate}
\usepackage[pagebackref]{hyperref} 

\newcommand{\FF}{{\mathbb{F}}}
\newcommand{\QQ}{{\mathbb{Q}}}
\newcommand{\ZZ}{{\mathbb{Z}}}

\newcommand{\bB}{{\mathbf{B}}}
\newcommand{\bG}{{\mathbf{G}}}

\newcommand{\cB}{{\mathcal{B}}}

\newcommand{\cO}{{\mathcal{O}}}

\newcommand{\fS}{{\mathfrak{S}}}

\newcommand{\ad}{{\operatorname{ad}}}
\newcommand{\End}{{\operatorname{End}}}
\newcommand{\Gal}{{\operatorname{Gal}}}
\newcommand{\Irr}{{\operatorname{Irr}}}
\newcommand{\SC}{{\operatorname{sc}}}
\newcommand{\Tr}{{\operatorname{Tr}}}
\newcommand{\Uch}{{\operatorname{Uch}}}

\newcommand{\SU}{{\operatorname{SU}}}
\newcommand{\SO}{{\operatorname{SO}}}
\newcommand{\Chevie}{{\sf{Chevie}}}

\newcommand{\wbG}{{\widehat{\bG}}}
\newcommand{\wG}{{\widehat{G}}}
\newcommand{\wH}{{\widehat{H}}}
\newcommand{\wL}{{\widehat{L}}}
\newcommand{\wM}{{\widehat{M}}}
\newcommand{\hrho}{{\widehat{\rho}}}

\newcommand{\tw}[1]{{}^{#1}\!}

\let\eps=\epsilon
\let\la=\lambda
\let\om=\omega
\let\si=\sigma

\newtheorem{thm}{Theorem}[section]
\newtheorem{lem}[thm]{Lemma}

\newtheorem{cor}[thm]{Corollary}
\newtheorem{prop}[thm]{Proposition}

\newtheorem{thmA}{Theorem}
\newtheorem{corA}[thmA]{Corollary}

\theoremstyle{definition}
\newtheorem{exmp}[thm]{Example}

\theoremstyle{remark}
\newtheorem{rem}[thm]{Remark}

\begin{document}

\title{Rationality of extended unipotent characters}

\author{Olivier Dudas}
\address{Institut de Mathématiques de Marseille -- I2M, Campus de Luminy,
         Avenue de Luminy, Case 930, 13288 Marseille Cedex 9, France}
\email{olivier.dudas@univ-amu.fr}

\author{Gunter Malle}
\address{FB Mathematik, RPTU Kaiserslautern, Postfach 3049,
         67653 Kaisers\-lautern, Germany.}
\email{malle@mathematik.uni-kl.de}

\begin{abstract}
We determine the rationality properties of unipotent characters of finite
reductive groups arising as fixed points of disconnected reductive groups under
a Frobenius map. In the proof we use realisations of characters in $\ell$-adic
cohomology groups of Deligne--Lusztig varieties as well as block theoretic
considerations.
\end{abstract}

\thanks{The second author gratefully acknowledges financial support by the DFG
-- Project-ID 286237555.}

\keywords{graph automorphism, rationality of extensions, unipotent characters}

\subjclass[2010]{20C08, 20C15, 20C33}

\date{\today}

\maketitle


\section{Introduction}   \label{sec:intro}
The work of George Lusztig has shown the singular importance of unipotent
characters in
the representation theory of finite reductive groups and hence of finite simple
groups of Lie type. Recent research has focussed on rationality properties of
characters of almost simple groups, and this naturally leads to the problem
of understanding fields of values, and fields of realisation, of extensions of
unipotent characters to groups of Lie type extended by graph or graph-field
automorphisms. The latter can be viewed as groups of rational points of suitable
disconnected algebraic groups. While the rationality properties of unipotent
characters themselves have long been known, due to the work of Lusztig
\cite{Lu02} and Geck \cite{Ge03} (see e.g.~\cite[Cor.~4.5.6]{GM20}), their
extensions to disconnected groups have so far not been studied systematically;
Digne--Michel \cite[Thm~II.3.3]{DM85} considered characters in the principal
series and \cite[Prop.~2]{Ma90} dealt with $\SU_3(q)$.

The field of values $\QQ(\rho)$ of a unipotent character $\rho$ is generated by
its Frobenius eigenvalue (see \cite[Prop.~4.5.5]{GM20}). Here, we show that this
Frobenius eigenvalue also governs the field of values of extensions of $\rho$.
Our first result concerns cuspidal characters:

\begin{thmA}   \label{thm:cusp}
 Let $\bG$ be a simple algebraic group with a Frobenius map $F$ and a commuting
 non-trivial graph automorphism $\si$. Then any cuspidal unipotent character
 $\rho$ of~$G=\bG^F$ has an extension $\hrho$ to $G\langle\si\rangle$ with
 $\QQ(\hrho)=\QQ(\rho)$, unless $G=\tw2A_{n-1}(q)$ with
 $n=\binom{t}{2}\equiv2,3\pmod4$ for some $t\ge3$, in which case
 $\QQ(\hrho)=\QQ(\sqrt{-q})$.\par
 In particular, $\QQ(\hrho)$ is generated by a $\delta$th root of the Frobenius
 eigenvalue of $\rho$, where $\delta\ge1$ is minimal such that $F^\delta$ acts
 trivially on the Weyl group of $\bG$.
\end{thmA}

Using earlier results on Frobenius--Schur indicators exhibits the following
connection:

\begin{corA}   \label{cor:FS}
 In the situation of Theorem~\ref{thm:cusp}, $\rho$ has a rational extension
 to~$G\langle\si\rangle$ if and only if $\rho$ has Frobenius--Schur
 indicator $+1$.
\end{corA}

For arbitrary unipotent characters, we obtain:

\begin{thmA}   \label{thm:HC}
 Let $\bG$ be a simple algebraic group with a Frobenius map $F$ and a commuting
 non-trivial graph automorphism $\si$. Then any $\si$-invariant rational
 unipotent character $\rho$ of~$G=\bG^F$ has a rational extension to
 $G\langle\si\rangle$, unless one of:
 \begin{enumerate}[\rm(1)]
  \item $G=A_{n-1}(q)$, $q$ is not a square and $\rho$ is labelled by a
   partition $\la=(\la_1,\ldots,\la_r)$ of~$n$ with
   \[\sum_i {{\la_i}\choose{2}}-\sum_i {{\la_i'}\choose{2}} + {{n}\choose{2}}
     \equiv1\pmod 2,\]
   where $\la'=(\la_1',\ldots,\la_s')$ is the partition conjugate to $\la$;
  \item $G=E_6(q)$, $q$ is not a square and $\rho$ is one of $\phi_{64,4}$
   or $\phi_{64,13}$; or
  \item $\rho$ lies in the Harish-Chandra series of a cuspidal unipotent
   character of a group of type $\tw2A_{n-1}(q)$ labelled by a $2$-core of size
   $n\equiv2,3\pmod4$.
 \end{enumerate}
 In cases~(1) and~(2), the extensions have character field $\QQ(\sqrt{q})$, in
 the third $\QQ(\sqrt{-q})$.
\end{thmA}

The case of cuspidal characters is settled in Section~\ref{sec:cusp}, where we
also give a reduction to simple groups, and in Section~\ref{sec:HC} we derive
the general case from the cuspidal one, thus proving
Theorem~\ref{thm:HC}. In Section~\ref{sec:exc} we discuss extensions of groups
of types $B_2$, $G_2$ and $F_4$ by an exceptional graph automorphism.

\section{Cuspidal unipotent characters}   \label{sec:cusp}

We consider the following setup. Let $\bG$ be a connected reductive linear
algebraic group with a Frobenius endomorphism $F:\bG\to\bG$ defining an
$\FF_q$-structure, and set $G:=\bG^F$, the finite group of $F$-fixed
points. We further assume that $\bG$ has a graph automorphism $\si$
commuting with~$F$. We set $\wbG=\bG\langle\si\rangle$ the semidirect
product of $\bG$ with~$\si$, and $\wG:=\wbG^F=G\langle\si\rangle$ the
corresponding extension of $G$.

\subsection{Deligne--Lusztig varieties and unipotent characters}
Let $\cB$ be the flag variety of $\bG$, that is, the variety of Borel subgroups
of $\bG$. The actions of $F$ and $\si$ on $\bG$ induce commuting
endomorphisms of $\cB$. The group $\bG$ acts by simultaneous conjugation on
$\cB\times\cB$ and the orbits are parametrized by the elements of the Weyl group
$W$ of $\bG$. Given $w \in W$, we denote by $\cO(w)$ the corresponding orbit and
we define the Deligne--Lusztig variety as in \cite[3.3]{Lu78} by
$$ X_w = \{ \bB\in\cB \mid (\bB,F(\bB)) \in \cO(w)\}.$$
The action of $\bG$ on $\cB\times\cB$ restricts to an action of the finite group
$G$ on $X_w$. Furthermore, $F$ (resp.~$\si$) induces a finite morphism
(resp.~an isomorphism) between $X_w$ and $X_{F(w)}$ (resp.~$X_{\si(w)}$).
Consequently, if $\delta$ is the smallest integer such that $F^\delta$ acts
trivially on $W$, then any Deligne--Lusztig variety has an action of $F^\delta$
commuting with $G$.

For $\ell$ a prime not dividing $q$ we denote by $R_w$ the corresponding
Deligne--Lusztig character of~$G$, given by
$$R_w(g)=\sum_{i\in\ZZ}(-1)^i\, \Tr\big(g \mid H^i_c(X_w,\QQ_\ell)\big)\qquad
  \text{for }g\in G.$$
This generalised character of $G$ does not depend on $\ell$ (see
e.g.~\cite[1.2]{Lu78}). The complex-valued characters which appear as
constituents of the various $R_w$ are called the \emph{unipotent characters}
of~$G$; we denote them $\Uch(G)$. By a result of Lusztig \cite[3.9]{Lu78}, for
any unipotent character $\rho$ of~$G$, $F^\delta$ acts by the same
\emph{eigenvalue of Frobenius} $\om_\rho$ on any $\rho$-isotypic component
$H^i_c(X_w,\overline{\QQ}_\ell)_\rho$ of any $\ell$-adic cohomology group of any
$X_w$, up to multiplication by integral powers of~$q^\delta$.

\subsection{Reduction to simple groups}
Here we follow the arguments in \cite[Rem.~4.2.1]{GM20} building upon
\cite[1.5.9--1.5.13]{GM20}. The centre $Z(\bG)$ is characteristic in~$\bG$, so
$\si$ induces a graph automorphism on $(\bG/Z(\bG))^F$. As unipotent characters
have $Z(\bG)^F$ in their kernel, for the purpose of studying rationality of
extensions of unipotent characters we may therefore assume $\bG$ is semisimple.
Let $\bG_\SC$ and $\bG_\ad$ be simply connected respectively adjoint groups of
the same type as $\bG$. Then these possess corresponding graph automorphisms
again denoted~$\si$. Furthermore, there are natural $F$- and $\si$-equivariant
epimorphisms $\bG_\SC\to\bG\to\bG_\ad$ such that the images of the respective
$F$-fixed points contain the derived subgroup of the $F$-fixed points of
the target. Now unipotent characters have the centre in their kernel and
restrict irreducibly to the derived subgroup. Since $\si$ stabilises the
centre and the derived subgroup, the rationality properties of extensions of
unipotent characters of $\bG^F$ agree with those of any group isogenous to it.
By passing to a group of adjoint type we see that we may hence assume for our
purposes that $\bG$ is simple (of a chosen isogeny type), which we do from now
on.

In particular, we can always assume that $\delta \in \{1,2,3\}$ and $F = \si^r$
in its action on $W$, for some $r\in\{1,2,3\}$ with one of $r$ or $\delta$ being
equal to~$1$.

\subsection{Eigenvalues of $F$ and character extensions}
We keep the above setting. We first look at the extensions over local fields
given by the $\ell$-adic cohomology of Deligne--Lusztig varieties.

\begin{prop}   \label{prop:Frob}
 Let $d$ be the order of $\si$ (recall that $d\in\{2,3\}$). Let $\rho\in\Uch(G)$
 be rational valued and $\si$-invariant. Assume that there is $w\in W$ such that
 \begin{enumerate}[\rm(1)]
  \item the $\langle F\rangle$-orbit of $w$ has length~$\delta$ and is
   $\si$-stable; and
  \item the multiplicity of $\rho$ in $R_w$ is not divisible by $d$.
 \end{enumerate}
 Then for every extension $\hrho$ of $\rho$ to $G\langle\si\rangle$, the field
 of values $\QQ_\ell(\hrho)$ is contained in the splitting field of
 $x^d-\om_\rho^{d/\delta}$ over $\QQ_\ell$. Furthermore, there is at least one
 extension $\hrho$ which is $\QQ_\ell$-valued if and only if there is a
 $\delta$th root of $\om_\rho$ in~$\QQ_\ell$.
\end{prop}

\begin{rem}
Observe that the conclusion of Proposition~\ref{prop:Frob} is well-defined
since the Frobenius eigenvalue $\om_\rho$ is uniquely determined up to integral
powers of $q^\delta$.
\end{rem}

\begin{proof}
We consider the subvariety
$$X = X_w \sqcup X_{F(w)} \sqcup \cdots \sqcup X_{F^{\delta-1}(w)}$$
of $\cB$. By~(1) it has an action of both $F$ and $\si$ and for all $i$ we have
$$H_c^i(X)_\rho \cong H_c^i(X_w)_\rho \oplus H_c^i(X_{F(w)})_\rho\oplus
  \cdots\oplus H_c^i(X_{F^{\delta-1}(w)})_\rho$$
as $\QQ_\ell G$-modules with $F$ cyclically permuting the $\delta$ summands.
By~(2) there is some~$i$ for which the multiplicity of $\rho$ in
$H_c^i(X_w)_\rho$ is not divisible by~$d$. Thus there also is a generalised
eigenspace $H_c^i(X_w)_{\rho,\mu}$ for an eigenvalue $\mu$ of $F^\delta$ on
$H_c^i(X_w)_\rho$ with the same property. Here, as cited above, $\mu$ differs
from $\om_\rho$ by an integral power of~$q^\delta$. 
First assume $\delta=1$. Then $H:=H_c^i(X_w)_{\rho,\mu}$ is a
$\QQ_\ell\wG$-module in which not all extensions of $\rho$ can occur with the
same multiplicity. Since $d\in\{2,3\}$, at least one of the extensions must then
be distinguished by its multiplicity and thus have values in $\QQ_\ell$. The
others are obtained by tensoring with linear characters of~$\wG/G$, which have
values in the splitting field of~$x^d-1$, so of $x^d-\om_\rho^d$.

Now assume $\delta=d$. Then $F$ has characteristic polynomial $(x^d-\mu)^m$ in
its action on $\bigoplus_{j=0}^{\delta-1} H_c^i(X_{F^j(w)})_{\rho,\mu}$, with
$m = \dim  H_c^i(X_w)_{\rho,\mu}$ and $m/\rho(1)$ not divisible by~$d$.
Let $K$ be the splitting field of $x^d-\mu$ over $\QQ_\ell$. Then
$\Gal(K/\QQ_\ell)$ permutes the generalised $F$-eigenspaces as it permutes the
eigenvalues, that is, as it acts on the roots of $x^d-\mu$. If $\QQ_\ell$
contains a zero of $x^d-\mu$, there is a $\QQ_\ell$-rational eigenspace and we
can argue as before to see that $\rho$ has a $\QQ_\ell$-rational extension.

If $x^\delta-\om_\rho$ has no zero in $\QQ_\ell$ then $K/\QQ_\ell$
is an extension of degree $\delta$ (recall that $\delta\leq 3$) and the $\delta$
different generalised eigenspaces of $F$ are Galois conjugate over $\QQ_\ell$.
Thus, the same holds for the $\delta$ different extensions of $\rho$.
\end{proof}

We now lift the rationality properties to $\QQ$.

\begin{lem}
 Let $\rho$ be a rational valued unipotent character of $G$. Assume that for
 all but finitely many $\ell$, $\rho$ has an extension $\hrho_\ell$ to $\wG$
 which takes values in $\QQ_\ell$. Then $\rho$ has an extension to~$\wG$ which
 takes values in~$\QQ$.
\end{lem}

\begin{proof}
We argue by contradiction. Assume all extensions of $\rho$ are defined over a
proper extension $K$ of $\QQ$, generated by a root of $f\in\QQ[x]$, say. Since
the sum of the extensions has values in~$\QQ$ and $\delta\in\{2,3\}$ this means
that $K$ (and hence $f$) has degree~$\delta$ over~$\QQ$. Note that $K/\QQ$ is
abelian. Thus by Dirichlet there are infinitely many primes~$\ell$ such that $f$
is also irreducible over $\QQ_\ell$, that is, its roots generate an extension of
degree~$\delta$. Hence for those $\ell$, $\rho$ does not have any
$\QQ_\ell$-rational extension, contradicting our assumption.
\end{proof}

\subsection{Extensions of cuspidal unipotent characters}
We keep the above setting. For the classification and properties of unipotent
characters we refer the reader to \cite[\S\S4.3--4.5]{GM20}. Recall that
cuspidal unipotent characters have values in $\QQ(\om_\rho)$.

\begin{prop}   \label{prop:untwisted}
 Let $\bG$ be simple of type $D_n$ or $E_6$, $\si$ a graph automorphism of $\bG$
 and $F$ a commuting Frobenius map with $\delta=1$. Then any cuspidal unipotent
 character $\rho$ of $G=\bG^F$ has an extension to $\wG=G\langle\si\rangle$
 defined over $\QQ(\rho)$.
\end{prop}

\begin{proof}
We consider the various cases according to Lusztig's classification of cuspidal
unipotent characters. By \cite[Thm~4.5.11]{GM20} they are all $\si$-invariant.
If $G$ is of type $D_n$ with $n=(2t)^2$ for some $t\ge1$ and $o(\sigma)=2$, then
the class of $W$ labelled by the bi-partition $(-;4t-1,4t-3,\ldots,1)$ contains
$\si$-stable elements (since the centraliser in $W$ of $\sigma$, of type
$B_{n-1}$, contains the class labelled $(-;4t-1,4t-3,\ldots,3)$), and by
\cite[Prop.~2.14]{Lu02}, the unique cuspidal unipotent character $\rho$ occurs
exactly once in the corresponding Deligne--Lusztig character. Thus
Proposition~\ref{prop:Frob} applies to show that $\rho$ has rational extensions.
If $G$ is of type $D_4$ with $o(\sigma)=3$, then there exists a $\si$-stable
element $w\in W$ in the class labelled by the bi-partition~$(-;31)$, and using
\Chevie\ \cite{MChev} the unique cuspidal unipotent character $\rho$ of $G$
appears with multiplicity~1 in~$R_w$. Thus Proposition~\ref{prop:Frob} applies
again.
\par
Finally, for $\bG$ of type $E_6$, there are two cuspidal characters
$E_6[\theta],E_6[\theta^2]$ with Frobenius eigenvalue a primitive third root of
unity $\theta$, respectively $\theta^2$.
Let $w\in W$ be in class $E_6$. Then $w$ can be chosen $\si$-stable. Again using
\cite{MChev}, both $E_6[\theta]$ and $E_6[\theta^2]$ occur in $R_w$ with
multiplicity~1. Let $H$ be a cohomology group of $X_w$ containing $\rho$ with
odd multiplicity, for $\ell$ a prime with $\ell\equiv1\pmod3$, so
$\sqrt{-3}\in\QQ_\ell$. Since $\si$ fixes $w$, it acts on $H$ and so the
$\QQ_\ell\wG$-module $H_\rho$
contains the two extensions $\rho_1,\rho_2$ of $\rho$ to $\wG$ with different
multiplicities. Thus, they must be $\QQ_\ell$-rational. Since this is true for
all $\ell\equiv1\pmod3$, the character field of $\rho_i$ is contained in
$\QQ(\theta)=\QQ(\rho)$, hence equal to $\QQ(\theta)$.
\end{proof}

We now turn to the twisted groups, where $\delta>1$.

\begin{prop}   \label{prop:SU}
 Let $G=\tw2A_{n-1}(q)$ where $n=\binom{t+1}{2}$ with $t\ge2$, and
 let $\rho$ be the cuspidal unipotent character of $G$. Let $\si$ be the
 graph-field automorphism of $G$ of order~$2$. Then the two extensions of $\rho$
 to  $\wG=G\langle\si\rangle$ are rational-valued if $\binom{n}{2}$ is even,
 and are algebraically conjugate over $\QQ(\sqrt{-q})$ if $\binom{n}{2}$ is odd.
\end{prop}

\begin{proof}
The cuspidal unipotent character of $\tw2A_{n-1}(q)$ is labelled by the
staircase partition $\la=(t,t-1,\ldots,1)$ (see e.g.~\cite[Prop.~4.3.6]{GM20}).
By \cite[Rem.~(a) after Thm~3.34]{Lu78} and \cite[Prop.~4.3.7]{GM20}) the
Frobenius eigenvalue of $\rho$ is $\om_\rho=(-q)^{\binom{n}{2}}$ (up to
multiplication by powers of~$q^2$). An application of the Murnaghan--Nakayama
rule shows that the irreducible character of $\fS_n$ labelled by $\la$ takes
value~$\pm1$ on elements $w\in W=\fS_n$ of cycle type $(2t-1,2t-5,\ldots)$.
Since the multiplicities of unipotent characters in Deligne--Lusztig characters
in type~$A_{n-1}$ are given by the character table of $W$ (see
\cite[Cor.~2.4.19]{GM20}), this means that $\rho$ has multiplicity~$\pm1$ in the
Deligne--Lusztig character~$R_w$. Also, no conjugate of $w$ is centralised
by~$\si$, hence neither by~$F$. Thus the assumptions of
Proposition~\ref{prop:Frob} are satisfied with $\delta=2$ and the conclusion
follows.
\end{proof}

\begin{rem}
Alternatively, by the arguments given for the principal series case in the
proof of Theorem~\ref{thm:HC} below, 
the conclusion of Proposition~\ref{prop:SU} would follow if the
full endomorphism algebra of the cohomology of $X_{w_0}$, where $w_0\in W$ is
the longest element, is indeed the Iwahori--Hecke algebra of type $A_{n-1}$ at
parameter $-q$ as speculated in \cite[3.10(b)]{Lu78} (see also the general
conjectures in \cite[1B]{BM93}).
\end{rem}

\begin{prop}   \label{prop:twisted}
 Let $\bG$ be simple of type $D_n$ or $E_6$, $\si$ a graph automorphism of $\bG$
 and $F$ a commuting Frobenius map with $\delta=o(\si)$. Then any cuspidal
 unipotent character $\rho$ of $G=\bG^F$ has an extension to
 $\wG=G\langle\si\rangle$ with field of values $\QQ(\rho)$.
\end{prop}

\begin{proof}
First assume that $G=\tw2D_n(q)$ where $n=(2t+1)^2$ with $t\ge1$, and let $\rho$
be the cuspidal unipotent character of $G$. Let $w\in W$ be in the class
$(-;4t+1,4t-1,\ldots,1)$. An application of Asai's formula
\cite[Thm~4.6.9]{GM20} shows that $\rho$ appears with multiplicity~$\pm1$
in~$R_w$ (see also \cite[2.19]{Lu02}). By \cite[Thm~4.11]{GM03} the
Frobenius eigenvalue of~$\rho$ is $\om_\rho=1$, up to multiplication by integral
powers of~$q^2$, so an application of Proposition~\ref{prop:Frob}
allows us to conclude.

Now assume $G=\tw3D_4(q)$, with $\si$ the graph-field automorphism of~$G$ of
order~$3$. If $\rho$ is the cuspidal unipotent character $\tw3D_4[-1]$ of $G$
then (using \Chevie \cite{MChev}), $\rho$ occurs with multiplicity~$1$ in
$X_w$ for $w$ of type $F_4$. The class of $w$ is not $F$-stable, so
Proposition~\ref{prop:Frob} applies. In this case the eigenvalue of $F^3$ for
$\rho$ equals $\om_\rho=-1$ (up to integral powers of~$q^3$) by
\cite[(7.3)]{Lu76}. 
Next, let $\rho$ be the cuspidal unipotent character $\tw3D_4[1]$, with
Frobenius eigenvalue~$\om_\rho=1$, by \cite[Rem.~4.9]{GM03}. It occurs with
multiplicity~$1$ in $X_w$ for $w$ of type $F_4(a_1)$, which can be chosen not
$\si$-invariant. Again Proposition~\ref{prop:Frob} applies.

Finally, let $G=\tw2E_6(q)$. The cuspidal unipotent character $\rho=\tw2E_6[1]$
with $\om_\rho=1$ (by \cite[Rem.~4.9]{GM03}) appears with multiplicity~1 in the
Deligne--Lusztig character $R_w$ for $w$ in class~$3A_2$. Choosing $w$ not
$\si$-stable, we conclude as in the previous case. Let now
$\rho$ be one of the cuspidal unipotent characters $\tw2E_6[\theta]$,
$\tw2E_6[\theta^2]$ of $G$ with $\om_\rho=\theta$, $\theta^2$
respectively \cite[(7.4)(e)]{Lu76}. For these the claim follows precisely as for
the non-rational cuspidal unipotent characters of $E_6(q)$ in the proof of
Proposition~\ref{prop:untwisted}.
\end{proof}

We are now ready to show our first main result:

\begin{proof}[Proof of Theorem~\ref{thm:cusp}]
All relevant cases are covered by Propositions~\ref{prop:untwisted},
\ref{prop:SU} and~\ref{prop:twisted}.
\end{proof}

\begin{proof}[Proof of Corollary~\ref{cor:FS}]
If $\rho$ is not real-valued the assertion holds trivially.  We now discuss the
real-valued cuspidal unipotent characters. All of them are rational by
\cite[Cor.~4.5.6]{GM20}. The Frobenius--Schur indicators of all these characters
of untwisted groups are $+1$ by Lusztig \cite[Thm~0.2]{Lu02}, and by
Proposition~\ref{prop:untwisted} they possess rational extensions.
\par
By Ohmori \cite{Oh96} the cuspidal unipotent character $\rho$ of the unitary
group $G=\SU_n(q)$ with $n=t(t+1)/2$ has Frobenius--Schur indicator
$(-1)^{\lfloor n/2\rfloor}$. Now $\lfloor n/2\rfloor$ is the $\FF_q$-rank
of~$G$, which in turn is congruent modulo~$2$ to the exponent~$i$ in the
Frobenius eigenvalue $\om_\rho=(-q)^i$ of $\rho$ by
\cite[Rem.~(a) after Thm~3.34]{Lu78}. The claim in this case thus follows from
Proposition~\ref{prop:SU}, observing that $\lfloor n/2\rfloor$ and
$\binom{n}{2}$ have the same parity.
\par
For the orthogonal group $G=\SO_{2n}^-(q)$ where $n=(2t+1)^2$, the
Frobenius--Schur indicator of the cuspidal unipotent character equals $+1$ by
\cite[1.13]{Lu02}, the Frobenius--Schur indicator of the cuspidal unipotent
character $\tw2E_6[1]$ of $\tw2E_6(q)$ is $+1$ by \cite[6.2]{Ge03}, and
similarly the indicators of the two cuspidal unipotent characters
of~$\tw3D_4(q)$ are also $+1$ by \cite[(7.6)]{Lu76} and \cite[6.2]{Ge03}. So for
the latter groups we may conclude by Proposition~\ref{prop:twisted}.
\end{proof}

\section{Harish-Chandra theory}   \label{sec:HC}
We now consider arbitrary unipotent characters.

\begin{proof}[Proof of Theorem~\ref{thm:HC}]
Let $\rho\in\Uch(G)$ be $\si$-invariant, so it has an extension $\hrho$ to
$\wG=G\langle\si\rangle$.  First assume $\rho$ lies in the principal series,
so it occurs as constituent in the permutation module $\QQ_\ell[G/B]$, whence
$\hrho$ occurs in $M:=\overline\QQ_\ell[\wG/B]$. Now as a
$\overline\QQ_\ell\wG\times\End_{\overline\QQ_\ell\wG}(M)$-bimodule, $M$ decomposes as the direct sum of
irreducible submodules $M_\phi$ indexed by $\phi\in\Irr(\End_{\overline\QQ_\ell\wG}(M))$
affording $\chi_\phi\otimes\phi$ for some $\chi_\phi\in\Irr(\wG)$ in the
principal series. Clearly, any Galois automorphism of $\QQ_\ell$ permutes the
$M_\phi$ as it permutes the characters~$\phi$, and hence it permutes the
$\chi_\phi$ in the same way. So the rationality statement for $\hrho$ follows
from the corresponding one for the extended Hecke algebra in
\cite[Thm~II.3.3]{DM85}. (See also the proof of \cite[Prop.~5.5]{Ge03}.)

Now assume we are not in that case. Let $L\le G$ be a
(split) Levi subgroup and $\la$ a cuspidal
unipotent character of $L$ such that $\rho$ lies in the Harish-Chandra series
of $\la$, so $\langle\rho,R_L^G(\la)\rangle\ne0$. Thus, $\rho$ corresponds to
a character $\phi$ of the relative Weyl group $W':=W_G(L,\la)$ (see
\cite[Thm~3.2.5]{GM20}). It is known that this relative Weyl group is of type
$A_2$, $B_n$, $G_2$ or $F_4$ in the cases we consider and that $L$ can also be
chosen to be $\si$-stable (see e.g.~\cite[Tab.~4.8]{GM20}). Moreover, as
$\la$ is the unique cuspidal unipotent character of $L$, it is also $\si$-stable.
Furthermore, $\si$ acts trivially on $W'$ except possibly if $W'$ has
type~$A_2$. By \cite[Thm~3.3 and~(3.6)]{GP} any such $\phi$ is of parabolic type,
that is, there is a parabolic subgroup $W_1'$ of $W'$ such that
$\langle 1_{W_1'}^{W'},\phi\rangle=1$. Note that we may assume $W_1'$ is proper
in~$W'$ if $W'\ne1$. This is clear if $\phi\ne 1_{W'}$, and $\phi=1_{W'}$ occurs
with multiplicity~1 in the permutation character on any parabolic subgroup.
Let $M\ge L$ be the Levi subgroup of~$G$ corresponding to $W_1'$ and $\chi$ be
the unipotent character of $M$ in the Harish-Chandra series $(L,\la)$
corresponding to~$1_{W_1'}$. Then by the comparison theorem
\cite[Thm~3.2.7]{GM20} this means that $\langle R_M^G(\chi),\rho\rangle=1$.
\par
First assume $\si$ has order~2 and let $\rho_1,\rho_2$ be the two extensions
of $\rho$ to~$\wG$. We claim that $\rho_1$ has the same rationality property as
some extension $\la_1$ of $\la$ to $\wL=L\langle\si\rangle$. If $W'=1$ then $L=G$
and the claim is obvious. If $W'\ne1$ then $W_1'<W'$ as argued above. 
Exclude for the moment the case that $W'$ has type $A_2$. Then $\si$ stabilises
$(W_1',\chi)$, so $\chi$ has two extensions
$\chi_1,\chi_2$ to $\wM=M\langle\si\rangle$. Then
$$ \langle R_\wM^\wG(\chi_1+\chi_2),\rho_1\rangle
  = \langle R_\wM^\wG(\mathrm{Ind}_M^\wM \chi),\rho_1\rangle
  =\langle \mathrm{Ind}_G^\wG R_M^G( \chi),\rho_1\rangle 
  =\langle R_M^G( \chi),\rho\rangle = 1$$
and thus (after possibly interchanging $\rho_1$ and $\rho_2$)
$\langle R_\wM^\wG(\chi_1),\rho_1\rangle=1$ and
$\langle R_\wM^\wG(\chi_1),\rho_2\rangle=0$.
Consequently $\rho_1$ has the same rationality properties as
$\chi_1$, which by induction has the same rationality properties as $\la_1$,
showing our claim.
Now, we obtain the desired conclusion from the rationality properties of
extensions of the cuspidal unipotent character~$\la$ discussed in
Section~\ref{sec:cusp}.
\par
Now assume that $W'$ has type $A_2$. Then by \cite[Tab.~4.8]{GM20} we have
$G=E_6(q)$ and $\la$ is the cuspidal unipotent character of $L$ of type $D_4$.
Of the three characters $\rho$ of $G$ in this Harish-Chandra series, two have
multiplicity one in $R_L^G(\la)$ and thus by our previous argument, possess
rational extensions to $\wG$. The Deligne--Lusztig character $R_w$ for $w$ in
the class $E_6$ (with characteristic polynomial $\Phi_3\Phi_{12}$) contains the
third character $D_4,\!r$ with multiplicity~$-1$. Since the class $E_6$ contains
elements centralised by $\si$, Proposition~\ref{prop:Frob} shows that
$D_4,\!r$ has a rational extension to~$\wG$.
\par
If $\si$ has order~3, then necessarily $G$ has type $D_4$ or $\tw3D_4$. Here,
only the cuspidal characters do not lie in the principal series, and for
those, the claim was shown above.
\end{proof}

\begin{exmp}
 Let $\rho$ be any of the two unipotent characters of $G=\tw2E_6(q)$ in the
 Harish-Chandra series of type $\tw2A_5$. Then the extensions of $\rho$ to
 $\wG$ have character field $\QQ(\sqrt{-q})$, by Proposition~\ref{prop:SU} in
 conjunction with Theorem~\ref{thm:HC}. These are Ennola-dual to the principal
 series characters $\phi_{64,4}$ and $\phi_{64,13}$ of $E_6(q)$, thus we see
 that the occurring irrationalities in Theorem~\ref{thm:HC} do obey the
 Ennola principle.
\end{exmp}

We can also understand completely the situation for groups of type $D_4$
extended by its full group of grah-automorphisms:

\begin{cor}
 Let $G=D_4(q)$ and $\Gamma\cong\fS_3$ its full group of graph automorphisms.
 Then all $\Gamma$-invariant unipotent characters of $G$ possess a rational
 extension to $\wG=G.\Gamma$.
\end{cor}

\begin{proof}
Let $\rho\in\Uch(G)$ be $\Gamma$-invariant. Let $\si\in\Gamma$ have order~3. By
Proposition~\ref{prop:untwisted}, $\rho$ has one rational extension $\rho_1$ to
$G\langle\si\rangle$ and two algebraically conjugate ones. Thus $\rho_1$ must
be $\Gamma$-invariant as well and further extends to two characters
$\hrho_1,\hrho_2$ of $\wG$. The restrictions of these to $G\langle\tau\rangle$,
where $\tau\in\Gamma$ has order~2, are the two extensions of $\rho$ to
$G\langle\tau\rangle$, so rational again by Proposition~\ref{prop:untwisted}.
But then $\hrho_1,\hrho_2$ must also be rational valued.
\end{proof}

\section{Exceptional graph automorphisms}   \label{sec:exc}
The groups $B_2(2^{2f+1})$, $G_2(3^{2f+1})$, $F_4(2^{2f+1})$ with $f\ge0$
possess exceptional outer graph automorphisms of order~2 not induced by an
automorphism of the ambient algebraic group (in particular one can not use the
results in Section~\ref{sec:cusp}). In the smaller two cases, the rationality
properties of extended unipotent characters were determined by Brunat:

\begin{prop}[Brunat]
 \begin{enumerate}[\rm(a)]
  \item For $G=B_2(2^{2f+1})$ all four invariant unipotent characters have
   rational extensions to the extension of $G$ by the exceptional graph
   automorphism.
  \item For $G=G_2(3^{2f+1})$, eight unipotent characters are invariant under
   the exceptional graph automorphism. Of these, the characters labelled
   $\phi_{1,0}, \phi_{2,1},\phi_{1,6}$ and $G_2[-1]$ have rational extensions,
   $\phi_{2,2}$ has an extension with character field $\QQ(\sqrt{3})$, while
   $G_2[1],G_2[\theta]$ and $G_2[\theta^2]$ have extensions with character field
   $\QQ(\sqrt{-3})$.
 \end{enumerate}
\end{prop}

\begin{proof}
The first statement follows by inspection of \cite[Tab.~6]{Br06}, the second
from \cite[Tab.~11]{Br07}. For the characters in the principal series this also
already follows from \cite[Tab.~IV]{GKP}, using the arguments in the proof of
Theorem~\ref{thm:HC}. (Note that the entries in the first rows of both tables
in loc.~cit. should correctly read $u^{l(w)}$, as confirmed by the authors.)
\end{proof}

The following observation will be used to handle the case of $F_4$; it would
also apply to prove rationality in some of the earlier cases.

\begin{lem}   \label{lem:block}
 Let $G\unlhd\wG$ with $\wG/G$ cyclic, $\ell$ a prime not
 dividing $|\wG:G|$ and $P$ a Sylow $\ell$-subgroup of $G$. Assume that
 $\wG=GC_{\wG}(P)$. Then any irreducible character in the
 principal $\ell$-block $B_0$ of $G$ has an extension to $\wG$ with the same
 character field.
\end{lem}

\begin{proof}
The assumptions imply that any $\rho\in\Irr(B_0)$ has a unique extension
$\hrho$ in the principal $\ell$-block $\hat B_0$ of $\wG$. Since
$\Irr(\hat B_0)$ can be defined in terms of orthogonality relations on
$\ell'$-elements, it is stable under all Galois automorphisms, which implies
$\QQ(\hrho)=\QQ(\rho)$, as claimed.
\end{proof}

\begin{prop}   \label{prop:F4}
 For $G=F_4(q)$ with $q=2^{2f+1}$ and $\si$ the exceptional graph automorphism,
 twenty-one unipotent characters are $\si$-invariant. Of these, the characters
 $$\phi_{1,0},\,\phi_{4,1},\,\phi_{9,2},\,\phi_{12,4},\,\phi_{6,6'},\,
   \phi_{6,6''},\,\phi_{4,8},\,\phi_{9,10},\,\phi_{4,13},\,\phi_{1,24},\,
   B_2,\!1,\,B_2,\!r,\,B_2,\!\eps,\,F_4^I[1],\,F_4^{II}[1]$$
 possess rational extensions to $\wG=G\langle\si\rangle$. The character
 $\phi_{16,5}$ has extensions with character field $\QQ(\sqrt{2})$, the cuspidal
 characters $F_4[\pm i]$ have extensions with character field $\QQ(i)$, the
 cuspidal characters $F_4[\theta],F_4[\theta^2]$ have extensions with character
 field $\QQ(\sqrt{-3})$, and $F_4[-1]$ has extensions with character field $\QQ(\sqrt{-2})$.
\end{prop}

\begin{proof}
Let $\bG$ be of type $F_4$ with a Steinberg endomorphism $F_0$ such that
$G=\bG^F$ for $F=F_0^2$. We may assume that $F_0$ induces $\si$ on $G$. For the
characters in the principal series the claim follows from the character table
of the extended Hecke algebra in \cite[Tab.~VI]{GKP} as in the proof of
Theorem~\ref{thm:HC}. For the characters $B_2,\!1$, $B_2,\!\eps$ we can also
argue as in the proof of Theorem~\ref{thm:HC} since they are parametrised in
their Harish-Chandra series by linear characters of the relative Weyl group, of
type $B_2$ (see \cite[Tab.~4.8]{GM20}), which thus have multiplicity~1 in its
regular character.
\par
Now let $w\in W$ be a $\si$-invariant element, and $\mu$ be an eigenvalue
of $F$ on $H = H_c^i(X_w)$ for some $i$. Let $\la = \sqrt{\mu}$
be a root of $\mu$. Then the action of $G$ on 
$\wH:=(H \otimes_{\QQ_\ell}\QQ_\ell(\lambda))_\mu$ extends to an action
of $\wG$ where $\sigma$ acts as $\la^{-1} F_0$ and the traces of all
$g\sigma$, for $g\in G$, lie in~$\la\QQ_\ell$. Thus, any Galois automorphism
sending $\la$ to $-\la$ will interchange the multiplicities in $\wH$ of the
two extensions of any irreducible character $\rho$ of $G$. In particular, if
$\rho$ has odd multiplicity in $H$, then its two extensions are only defined
over $\QQ_\ell(\la)$.
\par Now the eigenvalues of $F$ can be computed by evaluating
\cite[III, Prop. 1.2 with Thm~1.3]{DM85} with \Chevie\ \cite{MChev}.
Specifically, if $w\in W$ is a $\si$-invariant regular element of order~8,
then $F_4[-1]$ occurs with odd multiplicity in some $H_c^i(X_w)$
and with an eigenvalue of $F$ equal to $-q^3$ .Since $q$ is an odd
power of $2$ we have $\QQ_\ell(\sqrt{-q^3}) = \QQ_\ell(\sqrt{-2})$ and
the previous argument shows that $F_4[-1]$ has an extension in 
$\wH$ with character field $\QQ(\sqrt{-2})$.
\par
For the two characters $\rho=F_4[\pm i]$, we use the same element $w$. In that
case the eigenvalues of $F$ are $\pm iq^3$. As before, this means that both
characters possess extensions to $\wG$ with values in $\QQ(\sqrt{2i})$. Now note
that $\sqrt{2i}=\pm(1+i)$ lies in $\QQ(i)$, the character field of $\rho$.
\par
Next, the two cuspidal characters $F_4^I[1]$, $F_4^{II}[1]$ appear with
odd multiplicity in the $q^6$-eigenspace of $F$ on $H_c^i(X_w)$ for $w$
a $\si$-stable element in class $D_4(a_1)$. Arguing as before we see that these
characters possess rational extensions.
\par
For $\rho$ one of $F_4[\theta],F_4[\theta^2]$ or $B_2,\!r$ we use a block
theoretic argument. Let $\ell>2$ be a prime dividing $q^4-q^2+1$. Then
$G^{F_0}=\tw2F_4(2^{2f+1})$ contains a Sylow $\ell$-subgroup $P$ of~$G$, so
 $\wG=GC_\wG(P)$. Since $\rho$ lies in the principal $\ell$-block
(see e.g.~\cite[Thm~2.1(4)]{HL98}), Lemma~\ref{lem:block} shows
it has an extension to $\wG$ with the same character field.
\end{proof}

\begin{rem}
 Block theory also offers a way to see that the cuspidal unipotent character
 $F_4[-1]$ of $G=F_4(q)$ has non-real extensions to $\wG$: Let $\ell>2$ be a
 prime dividing $q^4+1$ and $P$ a Sylow $\ell$-subgroup of $G$. Then $P$ has
 larger automiser in $\wG$ than in $G$. Thus the Brauer tree of the principal
 $\ell$-block $\hat B_0$ of~$\wG$ is obtained by unfolding the Brauer tree of
 the principal $\ell$-block of~$G$ around the exceptional vertex. Since $F_4[-1]$
 and the trivial character of~$G$ lie on opposite sides of the exceptional
 vertex by \cite[Thm~2.1(3)]{HL98}  and the trivial character certainly has real
 extensions to~$\wG$, the extensions of $F_4[-1]$ cannot lie on the real stem of
 the Brauer tree of $\hat B_0$ and thus can't be
 real-valued. This does, however, not exhibit the precise character field. 
 \par
 Mutatis mutandis, this consideration also applies to the cuspidal unipotent
 characters of $\tw2A_2(q)$ and $\tw2A_5(q)$, and to $G_2[1]$ of $G_2(q)$ (with
 suitably chosen~$\ell$).
\end{rem}

\begin{rem}
 In the spirit of Ennola duality, the character fields of the extensions of
 $\phi_{2,2}$ and $G_2[1]$ in type $G_2$, as well as of $\phi_{16,5}$ and
 $F_4[-1]$ in type $F_4$ should be considered as being $\QQ(\sqrt{\pm q})$,
 respectively, since these pairs of characters are mutually Ennola dual.
\end{rem}



\begin{thebibliography}{13}

\bibitem{BM93}
{\sc M. Brou\'e, G. Malle}, Zyklotomische Heckealgebren. \emph{Ast\'erisque}
  No.~{\bf212} (1993), 119--189.

\bibitem{Br06}
{\sc O. Brunat}, The Shintani descents of Suzuki groups and their consequences.
  \emph{J. Algebra \bf303} (2006), 869--890. 

\bibitem{Br07}
{\sc O. Brunat}, On the extension of $G_2(3^{2n+1})$ by the exceptional graph
  automorphism. \emph{Osaka J. Math. \bf44} (2007), 973--1023.

\bibitem{DM85}
{\sc F. Digne, J. Michel}, Fonctions $L$ des vari\' et\' es de Deligne--Lusztig
  et descente de Shintani. \emph{M\'em. Soc. Math. France (N.S.)} No. {\bf20}
  (1985). 

\bibitem{Ge03}
{\sc M. Geck}, Character values, Schur indices and character sheaves.
  \emph{Represent. Theory \bf7} (2003), 19--55. 

\bibitem{GKP}
{\sc M. Geck, S. Kim, G. Pfeiffer}, Minimal length elements in twisted conjugacy
  classes of finite Coxeter groups. \emph{J. Algebra \bf229} (2000), 570--600. 

\bibitem{GM03}
{\sc M. Geck, G. Malle}, Fourier transforms and Frobenius eigenvalues for finite
  Coxeter groups. \emph{J. Algebra \bf260} (2003), 162--193.

\bibitem{GM20}
{\sc M. Geck, G. Malle}, \emph{The Character Theory of Finite Groups of Lie
  Type. A Guided Tour}. Cambridge University Press, Cambridge, 2020.

\bibitem{GP}
{\sc M. Geck, G. Pfeiffer}, \emph{Characters of Finite Coxeter Groups and
  Iwahori--Hecke Algebras.} The Clarendon Press, Oxford University Press,
  New York, 2000.

\bibitem{HL98}
{\sc G. Hiss, F. L\"ubeck}, The Brauer trees of the exceptional Chevalley groups
  of types $F_4$ and $\tw2E_6$. \emph{Arch. Math. (Basel) \bf70} (1998), 16--21. 

\bibitem{Lu76}
{\sc G. Lusztig}, Coxeter orbits and eigenspaces of Frobenius. \emph{Invent.
  Math. \bf38} (1976/77), 101--159.

\bibitem{Lu78}
{\sc G. Lusztig}, \emph{Representations of Finite Chevalley Groups}. CBMS
  Regional Conference Series in Mathematics, 39. American Mathematical
  Society, Providence, R.I., 1978.

\bibitem{Lu02}
{\sc G. Lusztig}, Rationality properties of unipotent representations. \emph{J.
  Algebra \bf258} (2002), 1--22.

\bibitem{Ma90}
{\sc G. Malle}, Some unitary groups as Galois groups over $\QQ$.
  \emph{J. Algebra \bf131} (1990), 476--482.

\bibitem{MChev}
{\sc J. Michel}, The development version of the CHEVIE package of GAP3.
 \emph{J. Algebra \bf435} (2015), 308--336. 

\bibitem{Oh96}
{\sc Z. Ohmori}, The Schur indices of the cuspidal unipotent characters of the
  finite unitary groups. \emph{Proc. Japan Acad. Ser. A Math. Sci. \bf72}
  (1996), 111--113.

\end{thebibliography}
\end{document}